\documentclass[11pt]{article}
\usepackage{amsfonts}
\usepackage{amsmath}
\usepackage{epsfig}

\title{The Optimal Paper Moebius Band}
\author{Richard Evan Schwartz \thanks{Supported by N.S.F. Grant
    DMS-2102802, a Simons Sabbatical Fellowship, and a Mercator Fellowship.}}

\newtheorem{theorem}{Theorem}[section]
\newtheorem{proposition}[theorem]{Proposition}
\newtheorem{lemma}[theorem]{Lemma}

\def\startproof{{\bf {\medskip}{\noindent}Proof: }}

\def\endproof{$\spadesuit$  \newline}

\def\R{\mbox{\boldmath{$R$}}}%

\begin{document}
\maketitle

\begin{abstract}
We prove that 
  a smooth embedded paper Moebius band
  must have aspect ratio greater than $\sqrt 3$.
  We also prove that any sequence of smooth
embedded paper Moebius bands whose aspect ratio
converges to $\sqrt 3$ must converge, up to
isometry, to the triangular Moebius band.
These results answer the mimimum aspect ratio
question discussed by W. Wunderlich in 1962 and
prove the more specific conjecture of
B. Halpern and C. Weaver from 1977.
    \end{abstract}

\section{Introduction}

\subsection{The Triangular Moebius Band}

To make a paper Moebius band you give a
$1 \times \lambda$ strip of
paper an odd number of
twists and then join the ends together.
For long strips this is easy and for short
strips it is difficult or impossible.  Let me first discuss
a beautiful example known as the
{\it triangular Moebius band\/}.
Figure 1a shows the triangular Moebius band. It
is based on a $1 \times \sqrt 3$ strip.

\begin{center}
\resizebox{!}{1in}{\includegraphics{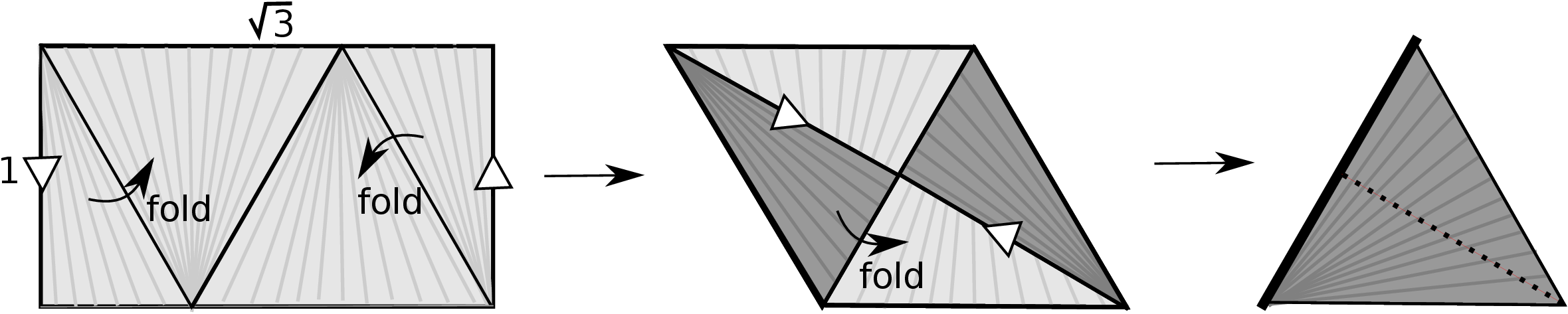}}
\newline
{\bf Figure 1a:\/} The triangular Moebius band
\end{center}

The strip in Figure 1a is lightly
shaded
on one side and darkly shaded
on the other.
First fold the flaps in to make a rhombus,
then fold the rhombus in half like a wallet.
This folding brings the two ends
together with a twist.  The dotted segment 
indicates where the ends are
joined. The bold segment
indicates the ``wallet fold''.  The dotted and bold segments
together make a pattern like a T.
The pinstriping exhibits the strip as a union of line segments,
disjoint except at the endpoints,
which stay straight during the folding.

\begin{center}
\resizebox{!}{1in}{\includegraphics{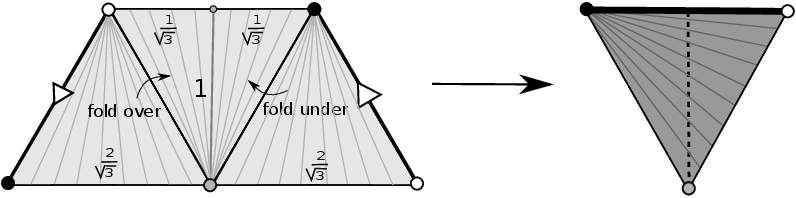}}
\newline
{\bf Figure 1b:\/} The triangular Moebius band: another view
\end{center}

Figure 1b shows another view.
Here we start
with a symmetric trapezoid rather than a rectangle,
but we get the same object when we fold
and join the sides together.   The bold edge indicates
where the sides are joined.
The dotted and bold segments again make a ``T-pattern''.

The triangular Moebius band
goes back at least to the 1930 paper
[{\bf Sa\/}] of M. Sadowski.   Technically, it does
not quite fit the definition of a (smooth, embedded)
paper Moebius band that we give below, but it is the limit of such.

\subsection{The Minimum Aspect Ratio Question}

The triangular Moebius band looks like an extremely efficient
construction.
Can we do better in terms of making $\lambda$ smaller?
To answer this question in a meaninful way, we first need a formal definition.
\newline
\newline
{\bf Definition:\/}
A {\it smooth paper Moebius band of
  aspect ratio $\lambda$\/} is an
infinitely differentiable isometric mapping
$I: M_{\lambda} \to \R^3$, where
$M_{\lambda}$ is the flat Mobius band
obtained by identifying the length-$1$ sides of a
$1 \times \lambda$ rectangle. That is:
\begin{equation}
  \label{quotient}
  M_{\lambda}=([0,\lambda] \times [0,1])\!/\sim, \hskip 30 pt
  (0,y) \sim (\lambda,1-y).
\end{equation}
An {\it isometric mapping\/} is a map which preserves arc-lengths.
The map is an {\it embedding\/} if it is injective,
and an {\it immersion\/} in general.
Let $\Omega=I(M_{\lambda})$.  We often write
$I: M_{\lambda} \to \Omega$.  We call $\Omega$
{\it embedded\/} when $I$ is an embedding.
\newline
\newline
\noindent
{\bf Remark:\/}
The smooth formalism  rules out examples
which render the main question meaningless.
For instance, you could fold any
rectangle (e.g. a square)
like an accordion
into a thin strip, twist, then tape.  This origami monster is not
the limit of smooth embedded paper Moebius bands.
We give an alternative definition of a paper
Moebius band a the top of \S \ref{bends0} which
is sufficient for our purposes and avoids smoothness.
\newline

The early papers of M. Sadowsky [{\bf Sa\/}] and
W. Wunderlich [{\bf W\/}] treat both the existence
and differential geometry of  paper Moebius bands.
(See [{\bf HF\/}] and [{\bf T\/}] respectively
for modern English translations.)
The paper [{\bf CF\/}] gives a modern differential
geometric framework for developable surfaces like $\Omega$.
The papers  [{\bf CK\/}],  [{\bf KU\/}], [{\bf RR\/}] and
[{\bf Sab\/}] are all studies of the differential geometry
of paper Moebius bands.
I learned about
paper Moebius bands from the great
expository article  [{\bf FT\/}, Chapter 14] by Dmitry Fuchs and
Sergei Tabachnikov.

W. Wunderlich discusses the minimum aspect ratio
question, without an explicit guess,
in the introduction of his 1962 paper [{\bf W\/}].
In their 1977 paper [{\bf HW\/}], Ben Halpern and Charles Weaver study
the minimum aspect ratio question in detail. They prove two things.
\begin{itemize}
\item For smooth immersed paper Moebius bands one has
  $\lambda>\pi/2$.    Moreover, for any $\epsilon>0$ one can find an
  immersed example with $\lambda=\pi/2+ \epsilon$.
\item There exists some $\epsilon_0>0$ such that
  $\lambda>\pi/2+\epsilon_0$ for a smooth
  embedded paper Moebius band.  This
  $\epsilon_0$ is not an explicit constant.
\end{itemize}
On the last line of
[{\bf HW\/}], Halpern and Weaver conjecture
that $\lambda>\sqrt 3$ for a smooth embedded
paper Moebius band.

\subsection{Results}
 \label{outline}

In this paper I will resolve the Halpern-Weaver conjecture.

\begin{theorem}[Main]
  A smooth embedded paper Moebius band has aspect ratio
  greater than $\sqrt 3$.
\end{theorem}

I will also prove that the triangular Moebius band is truly the optimal
paper Moebius band -- at least when it comes to aspect ratio
minimizing.

 \begin{theorem}[Triangular Limit]
   \label{converge}
   Let $\{I_n: M_{\lambda_n} \to \Omega_n\}$ be a sequence of smooth
   embedded paper
   Moebius bands with $\lambda_n \to \sqrt 3$.
   Then, up to isometry, $I_n$ converges
   uniformly to the map giving the triangular Moebius band.
 \end{theorem}

 \noindent
 {\bf Outline of the Proofs:\/}
 Let $I: M_{\lambda} \to \Omega$ be a smooth embedded paper Moebius band.
 A {\it bend\/} is
 line segment $B' \subset \Omega$ which cuts across $\Omega$
 and has its endpoints in the boundary.
 We call the pre-image
 $B=I^{-1}(B') \subset M_{\lambda}$ a {\it pre-bend\/}.   Since $I$ is
arc-length preserving and $B'$ is an arc-length minimizing
path between its endpoints, $B$ must also be an arc-length minimizing
path between its endpoints.  Hence $B$ is also a line segment.

 It is a classic
 fact that
 $\Omega$ has a continuous foliation $\beta$ by bends.
 See \S \ref{bends0}.
 The bends in  $\beta$ vary continuously
 and are pairwise disjoint.  The corresponding
 foliation of $M_{\lambda}$ by pre-bends is like
 the pinstriping in Figures 1a and 1b except that
 the pre-bends are disjoint even at the endpoints.
 
  We say that a
 $T$-{\it pattern\/} on $\Omega$
 is a pair of bends which lie in
 perpendicular intersecting lines.  Look again
 at the right sides of Figures 1a and 1b.
  We call the
  $T$-pattern {\it embedded\/} if the two bends are disjoint.
In \S \ref{lemmaTproof} we prove
\begin{lemma}[T]
  A smooth embedded paper Moebius band has
  an embedded $T$-pattern.
\end{lemma}
Here is the idea.
The space of pairs of unequal bends in $\beta$
has a $2$-point compactification which makes it into
the $2$-sphere, $S^2$.
We define a pair of odd functions on $S^2$
which detect a $T$-pattern when they have a common zero.
We apply the Borsuk-Ulam Theorem to get a
common zero.

In \S \ref{lemmaGproof} we prove
\begin{lemma}[G]
  A smooth embedded
  paper Moebius band with an embedded $T$-pattern has aspect
  ratio greater than $\sqrt 3$.
\end{lemma}
Here is the idea.
We choose an embedded $T$-pattern on $\Omega$ and
then cut $M_{\lambda}$ open along one of the
corresponding pre-bends.  The result
is a bilaterally symmetric trapezoid $\tau$.  See
Figure 2 below.  We
then solve an optimization problem which
involves mapping $\tau$ into $\R^3$ with
constraints coming from the geometry of trapezoids
and $T$-patterns.

The Main Theorem is an immediate consequence of
Lemma T and Lemma G.  
In \S \ref{seqmin} we prove the Triangular Limit Theorem by
examining what the proof of Lemma G
says about a minimizing sequence of examples.
\newline
\newline
{\bf Remarks:\/}
(1) If an embedded paper Moebius band has a long enough
strip that is contained in a single plane, it also has a
non-embedded $T$-pattern.  That is why we take special
care to speak of embedded $T$-patterns. \newline
(2)  The interested reader will be able to tweak our proof of
Lemma T to show that a smooth immersed
Moebius band has a $T$-pattern,  though not necessarily
an embedded one.  We do not want to fool around with this.
\newline
(3) Likewise, the interested reader would be able to tweak our
proof of Lemma G to prove that an immersed
paper Moebius band with an embedded $T$-pattern has
aspect ratio greater than $\sqrt 3$.  \newline
(4)
The ideas above are an outgrowth of my
earlier paper [{\bf S1\/}].   In [{\bf S1\/}] I prove a
version of Lemma T in a complicated way and with
the side hypothesis that $\lambda<7\pi/12$.
I then (correctly) deduce that
$\lambda \geq \phi= (1+\sqrt 5)/2$.
When I try to further improve this easy bound,
I  make an idiotic mistake:  I claim that when you
cut open $M_{\lambda}$ along a pre-bend you get a parallelogram.
This mistake  invalidates my final bound, a weird and forgettable
algebraic number
in $(\phi,\sqrt 3)$.
This paper supersedes [{\bf S1\/}]
and is independent from it.
\newline
(5)  My informal notes [{\bf S2\/}]
give  a slower and more elementary account of
my proofs. I designed [{\bf S2\/}] for
college students and advanced high school students
who want to learn the arguments.

\subsection{Additional Material and Context}

The proofs are done after \S 2, but I include some
more material in \S 3.

In \S \ref{remarks} I 
elaborate on some aspects of the proofs.

In \S \ref{path} I discuss an alternate framework for
Lemma T.

In \S \ref{fourpoint} I discuss some topics adjacent
to paper Moebius bands.  Let me also say a
few things here.
The paper [{\bf CKS\/}] and [{\bf DDS\/}] consider the related
question of tying a piece of rope into a knot
using as little rope as possible.   The papers
[{\bf D\/}] and [{\bf DL\/}] consider folded ribbon knots.
  [{\bf DL\/}, Corollary 25]
  is in some sense a special
  case of our two results, and
  [{\bf DL\/}, Conjecture 26] is a variant of the
  Halpern-Weaver Conjecture in the category of
  folded ribbon knots. Our Main Theorem
  incidentally resolves this folded ribbon knot conjecture.
  Some authors have considered ``optimal Moebius bands''
  from other perspectives, either algebraic
  [{\bf Sz\/}] or physical
  [{\bf MK\/}], [{\bf SH\/}].
  
  In \S \ref{multitwist}
  I discuss multi-twist Moebius bands and cylinders,
  and some new results about them which followed
  the writing of this paper, namely
  [{\bf BrS\/}], [{\bf H\/}] and  [{\bf S3\/}].

  \subsection{Acknowledgements}
  
I  thank Brienne Brown,
Matei Coiculescu,
Robert Connelly,
Dan Cristofaro-Gardiner, Elizabeth Denne,
Ben Halpern,
Dmitry Fuchs, Javier Gomez-Serrano, Aidan Hennessey,
Jan Neinhaus,
Anton Izosimov,
Jeremy Kahn, Rick Kenyon, Stephen D. Miller,
Noah Montgomery,
Sergei Tabachnikov, and
Charles Weaver for
helpful discussions about this subject.
I  especially
thank Matei for suggesting that I try for a
``mapping proof'' of Lemma T as opposed to the kind of
proof I had previously.  That suggestion led me to find
a really nice proof of Lemma T that greatly simplified
this paper.
Jeremy Kahn's remarks about the
Borsuk-Ulam Theorem also helped streamline the
proof of Lemma T.
Finally, I thank the anonymous referees for insightful and
 helpful comments which improved the exposition of the paper.

\newpage

    \section{Proofs of the Results}

  \subsection{Existence of a Bend Foliation}
  \label{bends0}

  Our proofs of the main results only use the fact that
  $I: M_{\lambda} \to \Omega$ is injective and arc-length preserving,
  and that $\Omega$ has a continuous foliation by bends.
The reader who prefers this alternate
definition can skip this section, in which we start with the
definition in the introduction and deduce the
bend foliation.

Let $U \subset \Omega$ denote the set of points with
nonzero mean curvature.  As with any smooth
surface-with-boundary the mean curvature makes
sense even at the boundary.

\begin{lemma}
  \label{mc}
  Each $p \in U$ lies in a unique bend $\gamma \subset U$.
\end{lemma}

\startproof
This is a classic result, a case of
either  [{\bf CL\/}, p. 314, Lemma 2] or
 [{\bf HN\/}, \S 3, Lemma 2].
 For simpler proofs see
 [{\bf Ma\/}] or my own notes  [{\bf S4\/}].

 Here is a quick sketch.
 Let $\eta$ denote the Gauss map. This is
 locally defined on $U$.
 Since $\Omega$ has zero Gauss curvature, and $U$ has nonzero
 mean curvature, the differential
 $d\eta$ has $1$-dimensional image and
 $1$-dimensional kernel at each tangent plane of $U$.
 Let $\gamma \subset U$ be the curve through $p$ and integral to ${\rm ker\/}(d\eta)$.
 Along $\gamma$,
 the triple $\{{\rm ker\/}(d\eta),{\rm im\/}(d\eta),\eta\}$
 defines an orthonormal frame.  $\eta$ is constant along $\gamma$ by
 integration.  Since $\gamma$ is part of a foliation by integral
 curves along which $\eta$ is constant,
 ${\rm im\/}(d\eta)$ is also constant along $\gamma$.
 Hence ${\rm ker\/}(d\eta)$, the tangent line to $\gamma$, is
 also constant along $\gamma$.  Hence $\gamma$ is a line segment.
  More work shows that $\gamma$ continues in either
 direction until reaching $\partial \Omega$.
   \endproof

\begin{proposition}
  \label{main0}
  $\Omega$ has a continuous foliation by bends.
\end{proposition}

\startproof
No two bends of $U$
can intersect, even at the boundary.
The reason: If they did intersect, then the
point of intersection would have
two distinct straight lines through it and
hence zero mean curvature.  This contradicts
Lemma \ref{mc}.  Hence $U$ has a
foliation by bends.  The disjointness of these bends
implies the continuity of the foliation.
  
Let $\tau$ be the closure of a component of
$\Omega-U$.
If $\tau$ has empty interior then $\tau$ is a line segment,
the limit of a sequence of bends. Hence
$\tau$ is also a bend.
    If $\tau$ has non-empty interior then $\tau$ is
    contained in a single plane because
         the Gauss map is constant on $\tau$.  Moreover,
         $\tau$ is a trapezoid: two opposite sides
         $\tau_1$ and $\tau_2$ are bends and the
         other two opposite sides lie in $\partial \Omega$.
      We foliate $\tau$ by bends in a canonical way, 
      interpolating continuously between $\tau_1$ and $\tau_2$:
      The interpolation
      uses bends which lie on either
    parallel  or coincident lines, depending on
    whether $\tau_1$ and $\tau_2$ lie on
    parallel or coincident lines.
    Doing this construction on all such components,
    we get our continuous foliation of $\Omega$.
    \endproof

  \subsection{Proof of Lemma T}
  \label{lemmaTproof}

  \noindent
  {\bf Definitions:\/}
  Let $I: M_{\lambda} \to \Omega$ be a smooth embedded
  paper Moebius band.  We choose a continuous
  foliation $\beta$ of $\Omega$ by
  bends, as guaranteed by Proposition \ref{main0}.
  The preimage $I^{-1}(\beta)$  is a continuous
  foliation of $M_{\lambda}$ by pre-bends. All
  our bends and pre-bends belong to these foliations.
  
Each bend $u$ has exactly $2$ unit
vectors $\pm \overrightarrow u$ parallel to it.
We call either one an {\it orientation\/} of $u$.
The {\it centerline\/} of $M_{\lambda}$ is the circle
$([0,\lambda] \times \{1/2\})/\!\sim$. 
 The {\it centerline
of $\Omega$\/} is the image of the centerline of
$M_{\lambda}$ under the map $I$.  These centerlines
are topological circles and $I$ is maps one to the other.
\newline
\newline
  \noindent
  {\bf Intersection with the Centerline:\/}
  Here is a proof that
  a pre-bend $u$ intersects the centerline of $M_{\lambda}$
    exactly once.  
  Let $\ell(\cdot)$ denote length.  If we have $\ell(u) < \sqrt{1+\lambda^2}$ we
can move $u$ by an isometry
so that it misses the vertical sides of $[0,\lambda] \times [0,1]$.
But then $u$ clearly intersects the centerline exactly once.
So, if $u$ intersects the centerline more than once, we have
$\ell(u) \geq \sqrt{1+\lambda^2}>\lambda$.
But $\partial \Omega=I(\partial M_{\lambda})$ is a loop that
contains the endpoints of the bend $u'=I(u)$.
Hence
$2\lambda=\ell(\partial \Omega) \geq 2\ell(u')=
2 \ell(u)>2\lambda$,
a contradiction.
\newline
\newline
\noindent
{\bf The Circle of Bends:\/}
    We parametrize
    the bends of $\beta$ by $\R/2\pi$, as follows:
Since $I$ is an embedding, our result about
pre-bends implies that
each bend of $\Omega$ intersects the centerline
of $\Omega$ exactly once.
    We associate to each bend the point where it intersects the
    centerline, and
  then we identify the centerline with $\R/2\pi$.
\newline
\newline
{\bf The Cylinder and the Sphere:\/}
Let $\Upsilon$ be the topological
cylinder of pairs $(x_0,x_1) \in (\R/2\pi)^2$ with $x_0 \not = x_1$.
A point $(x_0,x_1) \in \Upsilon$ corresponds to a pair $(u_0,u_1)$ of unequal
bends.
We let $\overline \Upsilon$ be the compactification of
$\Upsilon$ obtained by adding $2$ points:
$\partial_+$ (respectively
$\partial_-$) is the limit of pairs
$(x_0,x_1)$ where $x_1$ is just ahead (respectively just behind) $x_0$
in the cyclic order on $\R/2\pi$. The space
$\overline \Upsilon$ is homeomorphic to $S^2$, the $2$-sphere.
See \S \ref{remarks} for an explicit homeomorphism.
The map
$\Sigma(x_0,x_1)=(x_1,x_0)$ extends  to a continuous
involution of $S^2$ that swaps the two points $\partial_+$ and
$\partial_-$.
The explicit homeomorphism in \S \ref{remarks} conjugates
$\Sigma$ to the antipodal map of $S^2$.
\newline
\newline
{\bf Propagating the Orientations:\/}
    Let $(x_0,x_1) \in \Upsilon$. There is a unique path
    $t \to x_t$ in $\R/2\pi$ which joins $x_0$ to $x_1$, moves
    at constant speed, locally increases in the cyclic order
    on $\R/2\pi$, and has length less than $2\pi$.
    This path has length near $0$ (respectively near $2\pi$) when
$(x_0,x_1)$ is near
$\partial_+$ (respectively
$\partial_- $).
Let $u_t$ be the bend associated to $x_t$.
We write $\overrightarrow u_0 \leadsto \overrightarrow u_1$
when there is a continuous orientation of the bends $\{u_t\}$
that restricts to 
$\overrightarrow u_0$ and
$\overrightarrow u_1$.
Note that
$-\overrightarrow u_0\leadsto -\overrightarrow u_1$ and,
since $\Omega$ is a Moebius band, $\overrightarrow u_1 \leadsto
-\overrightarrow u_0$.
Also $\overrightarrow u_1$ converges to $\pm \overrightarrow u_0$ when
$(x_0,x_1)$ converges to $\partial_{\pm}$.
\newline
\newline
{\bf The Map:\/}
Let $m_j$ be the midpoint of $u_j$. Using
the dot product $(\cdot)$ and the cross product $(\times)$ define
$F=(g,h): \Upsilon \to \R^2$,  where
\begin{equation}
  \label{gh}
  g(x_0,x_1)=\overrightarrow u_0 \cdot \overrightarrow u_1,
  \hskip 30 pt
  h(x_0,x_1)=(m_0-m_1) \cdot (\overrightarrow u_0 \times \overrightarrow u_1).
\end{equation}
In this definition we mean that $\overrightarrow u_0 \leadsto \overrightarrow u_1$.
Our definition is independent of the chosen orientation since
$-\overrightarrow u_0 \leadsto -\overrightarrow u_1$.
Also, $F$ extends continuously
$S^2$ with $F(\partial_{\pm})=(\pm 1,0)$.
Since $\overrightarrow u_1 \leadsto -\overrightarrow u_0$ we have
$F \circ \Sigma=-F$.  Why?  Well, clearly $g(x_1,x_0)=-g(x_0,x_1)$, and
$$h(x_1,x_0)=(m_1-m_0) \cdot( \overrightarrow u_1 \times (-\overrightarrow u_0))=
      (m_1-m_0) \cdot (\overrightarrow u_0 \times \overrightarrow u_1)
      =-h(x_0,x_1).
      $$
      The Borsuk-Ulam Theorem says that $(0,0) \in F(S^2)$.
     Since
            $F(\partial_{\pm}) \not = (0,0)$ we have $(0,0) \in
            F(\Upsilon)$.
            \newline
            \newline
         {\bf Remark:\/} Here is a self-contained proof that
         $(0,0) \in F(\Upsilon)$.
         Suppose not.
      If $\gamma$ is a continuous path in $S^2$ which goes from
$\partial_+$ to $\partial_-$, then
$F(\gamma)$ goes from $(1,0)$ to $(-1,0)$, misses $(0,0)$, and
winds some half integer $w(\gamma)$ times around the origin.
All choices of $\gamma$ are homotopic to each other relative to
$\partial_{\pm}$, so
$w(\gamma)$ is independent of $\gamma$.
But consider $\gamma'=\Sigma(\gamma)$, re-oriented so that
it goes from $\partial_+$ to $\partial_-$.   Since
$F \circ \Sigma=-F$ the image $F(\gamma')$ is obtained
by rotating $F(\gamma)$ by $180$ degrees about $(0,0)$ then
re-orienting it so that it goes from $(1,0)$ to $(-1,0)$.
But then $w(\gamma')=-w(\gamma)$, a contradiction.
\newline
\newline
      {\bf Endgame:\/}
      Let $(u_0,u_1)$ be the bends corresponding to
      $(x_0,x_1) \in F^{-1}(0,0)$.
       First, $u_0$ and $u_1$ are disjoint because they belong to the same foliation.
        Second, $\overrightarrow u_0$ and $\overrightarrow u_1$ are orthogonal
        because $g(x_0,x_1)=0$.
        Third,  $\overrightarrow u_0$ and $\overrightarrow u_1$
        and
        $m_0-m_1$ are all orthogonal to $\overrightarrow u_0 \times \overrightarrow u_1$ because
        $h(x_0,x_1)=0$, and this easily implies that $u_0,u_1$ are
        coplanar.
        Hence $(u_0,u_1)$ is an embedded T-pattern.
         This proves    Lemma T.

\subsection{Proof of Lemma G}
\label{lemmaGproof}

Let $\ell$ denote arc length.
Let $\bigtriangledown$ be a triangle with horizontal base.
Let $\vee$ be the union of the two non-horizontal sides of
$\bigtriangledown$.

\begin{lemma}
  \label{A1}
If $\bigtriangledown$ has base $\sqrt{1+t^2}$ and height $h
\geq 1$ then
$\ell(\vee) \geq \sqrt{5+t^2}$, with equality iff
$\bigtriangledown$ is isosceles and $h=1$.
\end{lemma}

\startproof
This is an exercise in high school geometry.
Let $p_1,p_2,q$ be the vertices of $\bigtriangledown$,
with $p_1,p_2$ lying on the base.
Let $p_2'$ be the reflection of $p_2$ in the horizontal line through
$q$.
Note that $\bigtriangledown$ is isosceles iff $p_1,q,p_2'$ are
collinear.  By symmetry, the triangle inequality, and the Pythagorean
Theorem,
$$\ell(\vee)=\|p_1-q\|+\|q-p_2'\| \geq \|p_1-p_2'\|=\sqrt{1+t^2+4h^2}
\geq \sqrt{5+t^2}.$$
We get equality if and only if $p_1,q,p_2'$ are collinear and $h=1$.
\endproof

Let $I: M_{\lambda} \to \Omega$ be a smooth embedded paper Moebius band with an
embedded $T$-pattern.  Let
$S'=I(S)$ for any $S \subset M_{\lambda}$.
We have $\ell(\gamma)=\ell(\gamma ')$ for any curve
$\gamma \subset M_{\lambda}$.
We rotate $\Omega$ so that one of the bends of the
$T$-pattern, $T'$, lies in
$X$-axis and the other bend, $B'$, lies in the negative ray of the
$Y$-axis.  Next, we let
$B=I^{-1}(B')$ and $T=I^{-1}(T')$ be the corresponding pre-bends.

We cut $M_{\lambda}$ open along
$T$ to get a bilaterally
symmetric trapezoid $\tau$.   We normalize $\tau$
so that the parallel
sides  are horizontal.
Reflecting $\tau$ in the coordinate axes if needed,
we arrange that $u,v,w,x$ are
mapped to $\Omega$ as in Figure 2.
Compare Figure 1b.
 The quantities $t$ and $b$ (which are both positive in the case
 depicted)
 respectively denote the horizontal
 displacements of $T$ and $B$.

\begin{center}
\resizebox{!}{1.5in}{\includegraphics{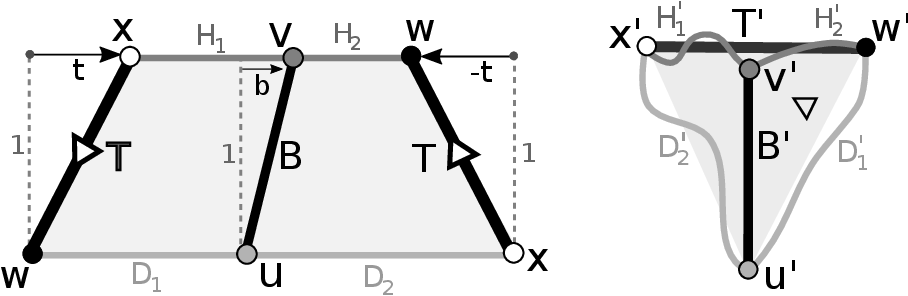}}
\newline
    {\bf Figure 2:\/} The trapezoid $\tau$ (left) and the T-pattern (right).
  \end{center}

  With Lemma \ref{A1} in mind, note that
  the shaded triangle $\bigtriangledown$ has
  $${\rm base\/}=\ell(T')=\ell(T)=\sqrt{1+t^2}, \hskip 20 pt
  {\rm height\/}>\ell(B')=\ell(B) =\sqrt{1+b^2} \geq 1.$$
  Let $H=H_1 \cup H_2$ and $D=D_1 \cup D_2$.
  We have $H', D' \subset \R^3$, so you should
  imagine you are floating in space above Figure 2 and looking down.
  Now $H'$ connects the endpoints of $T'$. Also,
  $D'$ connects the endpoints of $T'$ and contains $u$.
   From this structure, and from Figure 2 (left), we have
 \begin{align}
   \label{key01}
   \nonumber \ell(H)+\ell(D)=2\lambda. \\
   \nonumber
   \ell(D)-2t = \ell(H).\\
   \nonumber
   \sqrt{1+t^2}=\ell(T') \leq \ell(H')=\ell(H). \\
  \sqrt{5+t^2}< ^* \ell(\vee) \leq \ell(D')=\ell(D).
 \end{align}
The starred inequality is Lemma \ref{A1}.
Equation \ref{key01} give us
\begin{align}
  \label{con1}
  \nonumber
  \delta(t):=\sqrt{1+t^2}+\sqrt{5+t^2} <
  \ell(H)+\ell(D) = 2\lambda. \\
    \nu(t):=2 \sqrt{5+t^2}-2t < 2\ell(D)-2t =
  \ell(D)+\ell(H)=2\lambda.
\end{align}
Hence
\begin{equation}
  \label{max}
  2\lambda>\max(\delta(t),\nu(t)).
  \end{equation}
  Let $t_0=1/\sqrt 3$.    We have $\delta(t_0)=\nu(t_0)=2\sqrt 3$.
  Note the following.
  \begin{itemize}
    \item  $\delta$ is increasing on $(0,\infty)$.
      Hence $\delta(t)>2\sqrt 3$ if $t>t_0$.
    \item $\nu$ is decreasing on $\R$.  Hence
      $\nu(t)>2\sqrt 3$ if $t<t_0$.
  \end{itemize}
  Hence $\lambda > \sqrt 3$.  This proves Lemma G and, consequently,
  the Main Theorem.
    
\subsection{Proof of the Triangular Limit Theorem}
\label{seqmin}

Figure 3 shows what Figure 2 looks like with respect to
the triangular Moebius band and its T-pattern in Figure 1b.

\begin{center}
\resizebox{!}{1.3in}{\includegraphics{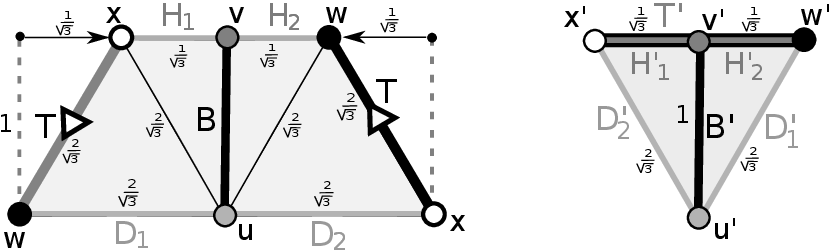}}
\newline
{\bf Figure 3:\/} Figure 2 for the triangular Moebius band.
\end{center}

Let $\{I_n: M_{\lambda_n} \to \Omega_n\}$ be as in the Triangular
Limit Theorem.   We run the construction of Lemma G for each
$n$ and analyze what happens as $n \to \infty$.
We use the same notation as in the proof of Lemma G, except
that we add subscripts to indicate the dependence on $n$.
\newline
\newline
\noindent
{\bf The Range:\/}
Since
$\max(\delta(t_n),\nu(t_n)) \to 2\sqrt 3$ we have
$t_n \to 1/\sqrt 3$.  Hence
$\ell(T'_n) \to 2/\sqrt 3$.
Since $\lambda_n \to \sqrt 3$ and
$\ell(\bigtriangledown_n) \leq 2\lambda_n$, we have
\begin{equation}
  \label{limsup}
  \limsup \ell(\bigtriangledown_n) \leq 2\sqrt 3.
  \end{equation}
The base of $\bigtriangledown_n$ converges to $2/\sqrt 3$ and
the height is always greater than $1$.
Lemma \ref{A1} combines with Equation \ref{limsup} to show that
$\bigtriangledown_n$ converges (up to isometries)
to an isosceles
triangle of base $2/\sqrt 3$ and height $1$.  But such
a triangle is actually equilateral.  This is
the shaded equilateral triangle $\bigtriangledown$
shown in Figure 3 (right).   We normalize by isometries so that we get actual convergence.
\newline
\newline
{\bf The Domain:\/}
Since $\ell(\partial \Omega_n)-\ell(\bigtriangledown_n) \to 0$ and also
$v_n'$ converges to the midpoint of $T'_n$, all the slack goes out of
Equation \ref{key01}, and
\begin{equation}
  \label{linear}
  \lim \ell(H_{n,1}')=
  \lim \ell(H_{n,2}')=1/\sqrt 3, \hskip 10 pt 
  \lim \ell(D_{n,1}')=
  \lim \ell(D_{n,2}')=2/\sqrt 3.
\end{equation}
Since $\ell(H_{n,1})=\ell(H_{n,1}')$, etc.
$\tau_n$ converges (up to isometries) to $\tau$, the trapezoid in
Figure 3 (left).
We normalize so that we get actual convergence.
\newline
\newline
{\bf The Boundary Map:\/}
The arcs $H_{n,1}', H_{n,2}', D_{n,1}',D_{n,2}'$ converge as sets
to the line segments connecting their endpoints because
$\ell(\partial \Omega_n)-\ell(\bigtriangledown_n) \to 0$.
Since $I_n$ is length preserving,
$I_n$ converges uniformly to a linear isometry when
restricted to each of $H_{n,1}, H_{n,2}, D_{n,1}, D_{n,2}$.
 Hence
 $I_n$ converges on $\partial M_{\lambda_n}$
 to the piecewise linear isometry
associated to the triangular Moebius band.
\newline
\newline
{\bf The Whole Map:\/}
We divide $M_{\lambda_n}$ into $3$ triangles. See Figure 4.
\begin{center}
\resizebox{!}{1.2in}{\includegraphics{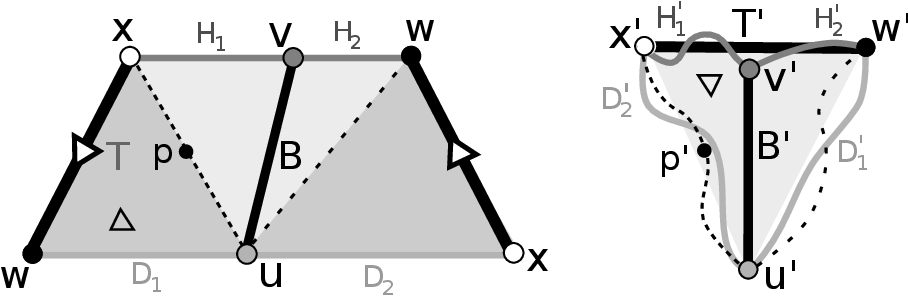}}
\newline
{\bf Figure 4:\/} Figure 2 revisited
  \end{center}
Consider the
restriction of $I_n$ to the left
triangle $\bigtriangleup_n$.   First,
$I_n$ is a linear isometry on
  $\overline{w_nx_n}=T_n$.
Second, $I_n$ converges to
a linear isometry on the segment
$\overline{w_nu_n}=D_{n,1}$.
Third,   $I_n$ converges to a linear isometry on
$\overline{x_nu_n}$ because $\|x_n'-u_n'\| \to \|x_n-u_n\|$ and
$I_n$ is distance non-increasing.
In summary, $I_n$ converges to a linear isometry on $\partial
\bigtriangleup_n$.
Since $I_n$ is distance non-increasing, this implies that
$I_n$ converges to a linear isometry on $\bigtriangleup_n$.
The same argument works for the other $2$ triangles.
Hence $I_n$ converges to the piecewise linear isometry associated to the
triangular Moebius band.
\newline
\newline
{\bf The Foliations:\/}  Our proof of the Triangular Limit Theorem is
done, but there is more to say.
Let $p_n$ be the midpoint of $\overline{x_nu_n}$.  The point
$p_n'$ converges to the midpoint of $\overline{x_n'u_n'}$.
The bend $\gamma_n'$ through $p_n'$ has its endpoints on
$H_{n,1}'$ and $D_{n,1}'$, and
this forces one endpoint of $\gamma_n'$ to converge to $x_n'$ and the other to
$u_n'$.  This means the pre-bend $\gamma_n=I_n^{-1}(\gamma_n')$ converges
to $\overline{x_nu_n}$.   In other words, the left dotted line in
Figure 4 is, for large $n$, quite close to a pre-bend.
The same argument works for the right dotted line.  From this, we
see that the pre-bend foliation of $M_{\lambda_n}$
converges to the pinstriping shown in Figures 1a (left) and 1b (left).
Likewise, the bend foliation of $\Omega_n$ converges to
the pinstriping shown in Figures 1a (right) and 1b (right).

\newpage

    \section{Discussion}
\label{discuss}

\subsection{An Explicit Homeomorphism}
\label{remarks}

Referring to the proof of Lemma T,
here we give an explicit homeomorphism from
$S^2$ to $\overline \Upsilon$.
Let $X_{\theta}$ denote the bend corresponding to
$\theta \in \R/2\pi$.

We
identify $\partial_+$ and $\partial_-$ respectively with the north and south pole of
$S_2$.  We parametrize $S^2-\partial_{\pm}$ by
    $(\theta,\phi)$,  where $\theta \in \R/2\pi$ is the
    longitude and $\phi \in (0,\pi)$, the angle with the vector
    pointing to
    $\partial_+$, is the latitude.  The antipodal map is
    $(\theta,\phi) \to (\pi+\theta,\pi-\phi)$ in these coordinates.
    We get a homeomorphism between
    $S^2-\partial_{\pm}$ and $\Upsilon$ with the correspondence
      $$(\theta,\phi) \leftrightarrow
      (X_{\theta-\phi},X_{\theta+\phi}).$$
      This conjugates
      $\Sigma$ to the antipodal map, and extends to a
      homeomorphism between $S^2$ and $\overline \Upsilon$.

\subsection{Paths of Oriented Lines}
\label{path}

    Anton Izosimov and Sergei Tabachnikov
    independently suggested to me the following
    generalization of Lemma T.
    \begin{lemma}
      \label{TT}
  Suppose $\{L_t|\ t \in [0,1]\}$ is a continuous family of
  oriented lines in $\R^3$ such that $L_1=L_0^{\rm opp\/}$, the same
  line as $L_0$ but with the opposite orientation. Then
  there exist parameters $r,s \in [0,1]$ such that
  $L_r$ and $L_s$ are perpendicular intersecting lines.
\end{lemma}

\startproof
This has the same proof as Lemma T, once we
  observe that our function $h$, defined in
  Equation \ref{gh}, is more natural than
  we have let on.  The points $m_0,m_1$ in the definition
  of $h$ could be any points on
  $u_0,u_1$ and we would get the same result.
  $g$ and $h$ are invariants of pairs of oriented lines.
  \endproof
  
Sergei also suggested to me a beautiful alternate formalism
for Lemma T: the {\it dual numbers\/}.
These have the form
$x + \epsilon y$ where $x,y \in \R$ and $\epsilon^2=0$.
Relatedly, the
 {\it dual vectors\/}
have the form $\overrightarrow a + \epsilon \overrightarrow b$,
where $\overrightarrow a, \overrightarrow b \in \R^3$ 
and again $\epsilon^2=0$.  In this context, the
dot product of two dual vectors makes sense
as a dual number.  See [{\bf HH\/}] for an exposition.

Each oriented line $\ell \subset \R^3$ gives rise to a
dual vector $\xi_{\ell}=\overrightarrow a + \epsilon \overrightarrow b$
where $\overrightarrow a$ is the unit vector pointing in the
direction of $\ell$ and $\overrightarrow b=\ell' \times \overrightarrow a$.
Here $\ell' \in \ell$ is any point.  All choices of $\ell'$ give rise
to the same $\overrightarrow b$; this vector is called the
{\it moment vector\/} of $\ell$.  This formalism identifies
the space of oriented lines in $\R^3$ with the so-called
{\it Study sphere\/} consisting of dual vectors
$\xi$ such that $\xi \cdot \xi = 1$.
The dual dot product $\xi_{\ell} \cdot \xi_{m}$ vanishes
if and only if $\ell$ and $m$ are perpendicular and intersect.

\subsection{Related Topics}
\label{fourpoint}
\label{foldedribbon}

\noindent
{\bf Square Peg:\/}
Around the time I got interested in the
Halpern-Weaver Conjecture I had been thinking
about the Toeplitz Square Peg Conjecture, which asks if
every continuous loop in the plane contains $4$ points which make
the vertices of a square.  See [{\bf Mat\/}] for a fairly
recent survey.  One can view a
$T$-pattern as a collection of $4$ points in the boundary of
the Moebius band which satisfy certain additional constraints -- e.g.
they are coplanar.
Put this way, a $T$-pattern is sort of like a square inscribed in a
Jordan loop.
\newline
\newline
{\bf Quadrisecants:\/}
The idea for Lemma T is also
similar in spirit for the idea developed in [{\bf DDS\/}] concerning
$4$ collinear points on a knotted loop.  These so-called
{\it quadrisecants\/} play a role similar to Lemma T in getting
a lower bound for the length of a knotted rope.
\newline
\newline
{\bf Folded Ribbon Knots:\/}
Elizabeth Denne pointed out to me the connection between
paper Moebius bands and {\it folded ribbon knots\/}.  Her paper
with Troy Larsen [{\bf DL\/}] gives a formal definition
of a folded ribbon knot and has a wealth of
interesting constructions, results, and conjectures.
See also [{\bf D\/}], a survey article.

Informally, folded ribbon knots are the objects
you get when you take a flat cylinder or Moebius band,
fold it into a knot,
and then press it into the plane.  Associated to a
folded ribbon knot is a polygon, which comes from the
 centerline of the object.
Even though
the ribbon knot lies entirely in the plane, one assigns
additional combinatorial data which keeps track of
``infinitesimal'' under and over crossings as in a
knot diagram.  So the associated centerline is really a
knot (or possibly the unknot).

[{\bf DL\/}, Corollary 25] proves our Main Theorem in the
category folded ribbon Moebius bands whose associated
centerline is a triangle.  This is
a finite dimensional problem.  [{\bf DL\/}, Conjecture 26]
says that  [{\bf DL\/}, Corollary 25] is true without
the very strong triangle restriction.  This
is an infinite dimensional problem like
the Halpern-Weaver Conjecture.  
The combination of our Main Theorem and the Triangular Limit
Theorem implies   [{\bf DL\/}, Conjecture 26].  One takes 
arbitrarily nearby smooth approximations, as in [{\bf HW\/}],
and then applies our results to them.

\subsection{More Twists}
\label{multitwist}

The Halpern-Weaver Conjecture is one of infinitely many
similar kinds of questions one can ask about paper
Moebius bands.  E.g., one can translate
the many conjectures made in [{\bf DL\/}] about folded ribbon knots into
the broader language of paper geometry.
Here I will discuss some recent developments that followed
(the preprint version of) this paper.

    One can make a {\it twisted cylinder\/} by taking a $1 \times \lambda$
    strip of paper, giving it an even and nonzero number of twists, and
    then joining the ends together.
    The essential feature of twisted cylinders is that their two
    boundary components make a nontrivial link.
    There are two optimal limiting shapes which have
    interpretations as folded ribbon knots.   Both are
    folding patterns which wrap a $1 \times 2$ strip $4$ times around
    a right-angled isosceles triangle.  In [{\bf S3\/}] I prove that
    a twisted cylinder has aspect ratio greater than $2$ and that any
    minimizing sequence converges on a subsequence to one of the two
    optimal models.  This result also confirms the $n=1$ case of
    [{\bf DL\/}, Conjecture 39].
    The proof is somewhat similar to what I do in this paper, though
    the fine-scale details are different.  Noah Montgomery (private
    communication)
    independently
    came up with a proof of the cylinder result.  His elegant  proof is different than
    mine.

    An essential feature of $3$-twist paper Moebius bands
        is that their boundaries are knotted.
        Brienne Brown did some experiments with these objects and found
        two candidate optimal models.  We
         call these the {\it crisscross\/} and
        the {\it cup\/}.  Both are made from a $1 \times 3$ strip of
        paper. The crisscross is planar, and has an interpretation
        as a folded ribbon knot.  The cup is not-planar: It is a double
        wrap of $3$ mutually orthogonal right-angled isoceles
        triangles arranged like $3$ faces of a tetrahedron.  We
        wrote about this in [{\bf BrS\/}], and conjecture there that
        $\lambda>3$ for an embedded
        multi-twisted paper Moebius band.

        More recently, Aidan Hennessey [{\bf H\/}] proved
        the fantastic result that
            one can make a cylinder or a Moebius band with any
            number of twists using a $1 \times 6$ strip.
            Extremely recently, following my
            lecture at UCLA on 8 Oct 2024, I told  Jan Neinhaus about
            Hennessey's construction.  The next day,
            Jan showed me how one can optimize
            and get  $3\sqrt 3+ \epsilon$,
            for any $\epsilon>0$, in place of $6$.
            Jan's limiting shape is half a regular hexagon.  We
             conjecture that $3 \sqrt 3$ is optimal in this context.
             Nobody has written about this yet.
                         
\newpage

\section{References}

\noindent
    [{\bf BrS\/}] B. E. Brown and R. E. Schwartz,
    {\it The crisscross and the cup: Two short $3$-twist paper Moebius bands\/},
    preprint 2023, arXiv:2310.10000
\vskip 6 pt
\noindent
[{\bf CF\/}] Y. Chen and E. Fried,
{\it Mobius bands, unstretchable material sheets and developable surfaces\/},
Proceedings of the Royal Society A, (2016)
\vskip 6 pt
\noindent
[{\bf CK\/}] Carmen Chicone, N. J. Kalton, {\it Flat Embeddings of the
Moebius Strip in $\R^3$\/}, Communications in Nonlinear Analysis 9
(2002), no. 2, pp 31-50
\vskip 6 pt
\noindent
[{\bf CKS\/}] J. Cantarella, R. Kusner, J. Sullivan, {\it On the
minimum ropelength of knots and links\/},  Invent. Math. {\bf 150\/} (2) pp 257-286 (2003)
\vskip 6 pt
\noindent
  [{\bf CL\/}], S.-S. Chern and R. K. Lashof,
    {\it On the total curvature of immersed manifolds\/},
    Amer. J. Math. {\bf 79\/} (1957) pp 306--318
    \vskip 6 pt
    \noindent
        [{\bf D\/}] E. Denne, {\it Folded Ribbon Knots in the Plane\/}, 
    The Encyclopedia of Knot Theory (ed. Colin Adams, Erica Flapan, Allison Henrich, Louis H. Kauffman, Lewis D. Ludwig, Sam Nelson)
    Chapter 88, CRC Press (2021)
    \vskip 6 pt
    \noindent
    [{\bf DL\/}] E. Denne, T. Larsen, {\it Linking number and folded ribbon unknots\/},
    Journal of Knot Theory and Its Ramifications, Vol. 32 No. 1 (2023)
    \vskip 6 pt
    \noindent
        [{\bf DDS\/}] E. Denne, Y. Diao, J. M. Sullivan, {\it Quaadrisecants give new lower bounds for the ropelength of a knot\/},
        Geometry$\&$Topology 19 (2006) pp 1--26
        \vskip 6 pt
        \noindent
[{\bf FT\/}], D. Fuchs, S. Tabachnikov, {\it Mathematical Omnibus: Thirty Lectures on Classic Mathematics\/}, AMS 2007
\vskip 6 pt
\noindent
[{\bf H\/}] A. Hennessey, {\it Constructing many-twist M\"obius bands
  with small aspect ratios\/}.  arXiv:2401:14639, to appear in
Comptes Rendus
\vskip 6 pt
\noindent
[{\bf HF\/}], D.F. Hinz, E. Fried, {\it Translation of Michael Sadowsky’s paper 
‘An elementary proof for the existence of a developable MÖBIUS band and the attribution 
of the geometric problem to a variational problem’\/}. J. Elast. 119, 3–6 (2015)
\vskip 6 pt
\noindent
[{\bf HH\/}], A. A. Harkin and J. B. Harkin, {\it Geometry of
Generalized Complex Numbers\/}, Mathematics Magazine {\bf 77\/} (2),
pp 118-129  (2004)
\vskip 6 pt
\noindent
 [{\bf HL\/}], P. Hartman and L. Nirenberg,
  {\it On spherical maps whose Jacobians do not change sign\/},
    Amer. J. Math. {\bf 81\/} (1959) pp 901--920
    \vskip 6 pt
    \noindent
[{\bf HW\/}], B. Halpern and C. Weaver,
{\it Inverting a cylinder through isometric immersions and embeddings\/},
Trans. Am. Math. Soc {\bf 230\/}, pp 41--70 (1977)
\vskip 6 pt
\noindent
[{\bf KU\/}], Kurono, Yasuhiro and Umehara, Masaaki,
{\it Flat M\"obius strips of given isotopy type in $\R^3$ whose
  centerlines are
  geodesics or lines of curvature\/}, Geom. Dedicata 134 (2008), pp
109-130
\vskip 6 pt
\noindent
[{\bf MK\/}] L. Mahadevan and J. B. Keller, 
{\it The shape of a M\"obius band\/}, Proceedings of the Royal Society A (1993) 
\vskip 6 pt
\noindent
[{\bf Mas\/}] W. S. Massey, {\it Surfaces of Gaussian Curvature Zero
  in Euclidean $3$-Space\/},  Tohoku Math J. (2) 14 (1), pp 73-79
(1962)
\vskip 6 pt
\noindent
[{\bf Ma1\/}] B. Matschke, {\it A survey on the
Square Peg Problem\/},  Notices of the A.M.S.
{\bf Vol 61.4\/}, April 2014, pp 346-351.
\vskip 6 pt
\noindent
[{\bf RR\/}] T. Randrup, P. Rogan, {\it Sides of the M\"obius
  Strip\/},
Arch. Math. 66 (1996) pp 511-521
\vskip 6 pt
\noindent
[{\bf Sa\/}], M. Sadowski, {\it Ein elementarer Beweis für die Existenz eines 
abwickelbaren MÖBIUSschen Bandes und die Zurückführung des geometrischen 
 Problems auf einVariationsproblem\/}. Sitzungsberichte der Preussischen Akad. der Wissenschaften, physikalisch-mathematische Klasse 22, 412–415.2 (1930)
\vskip 6 pt
\noindent
[{\bf Sab\/}] I. Kh. Sabitov, {\it Isometric immersions and embeddings
  of a flat M\"obius strip into Euclidean spaces\/}, Izv. Math. {\bf
  71\/} (2007), pp 1049-1078
\vskip 6 pt
\noindent
[{\bf S1\/}] R. E. Schwartz, {\it An improved bound on the optimal paper Moebius band\/},
Geometriae Dedicata, 2021
\vskip 6 pt
\noindent
[{\bf S2\/}] R. E. Schwartz, {\it The Optimal Paper Moebius Band: A
  Friendly Account\/}, preprint, 2023
\vskip 6 pt
\noindent
  [{\bf S3\/}] R. E. Schwartz, {\it The Optimal Twisted Paper Cylinder\/},
 preprint 2023, arXiv:2309.14033
 \vskip 6 pt
\noindent
[{\bf S4\/}] R. E. Schwartz, {\it Existence of Bends on Paper Moebius bands\/}
\newline
http://www.math.brown.edu/$\sim$ res:  Math Notes section
\vskip 6 pt
\noindent
[{\bf Sz\/}] G. Schwarz, {\it A pretender to the title ``canonical Moebius strip''\/},
Pacific J. of Math., {\bf 143\/} (1) pp. 195-200, (1990)
\vskip 6 pt
\noindent
    [{\bf SH\/}] E. L. Starostin, G. H. M. van der Heijden, {\it The shape of a M\"obius Strip\/},
    Nature Materials {\bf 6\/} (2007) pp 563 -- 567
    \vskip 6 pt
    \noindent
[{\bf T\/}] Todres, R. E., {\it Translation of W. Wunderlich's On a Developable M\"obius band\/},
Journal of Elasticity {\bf 119\/} pp 23--34 (2015)
\vskip 6 pt
\noindent [{\bf W\/}] W. Wunderlich, {\it \"Uber ein abwickelbares M\"obiusband\/}, Monatshefte f\"ur Mathematik {\bf 66\/} pp 276--289 (1962)

\end{document}